\documentclass[12pt,a4paper]{amsart}

\usepackage{amsmath}

\usepackage[english]{babel}

\def\AmS{{\protect\usefont{OMS}{cmsy}{m}{n}%
  A\kern-.1667em\lower.5ex\hbox{M}\kern-.125emS}}
\def\latex/{{\protect\LaTeX}}
\def\amslatex/{{\protect\AmS-\protect\LaTeX}}

\DeclareMathOperator{\ad}{ad}
\DeclareMathOperator{\AFS}{AFS}

\numberwithin{equation}{section}

\renewcommand{\theta}{\vartheta}

\newcommand{\p}{\text{\emph{p}\/}\nobreakdash}
\newcommand{\F}{\text{$\mathbf{F}$}}
\newcommand{\U}[2]{\underbrace{#1 \dots #1}_{#2}}
\newcommand{\Meta}{\mathcal{M}}
\newcommand{\Infl}[2]{{}^{#1} \! #2}
\newcommand{\uptoasign}{\approx}
\newcommand{\trivfrac}[2]{#1 / #2}

\newcommand{\From }{From } 
\newcommand{\SOTC}{sequence of two-step centralizers}
\newcommand{\TSC}{two-step centralizer}
\newcommand{\TSCs}{two-step centralizers}

\newcommand{\SOCL}{sequence of constituent lengths}

\newcommand{\ignore}[1]{}

\theoremstyle{plain}
\newtheorem{step}{Step}
\swapnumbers
\newtheorem{dummy}{}[section]
\newtheorem{prop}[dummy]{Proposition}
\newtheorem{theorem}[dummy]{Theorem}

\newtheorem{lemma}[dummy]{Lemma}

\theoremstyle{definition}
\newtheorem{definition}[dummy]{Definition}
\newtheorem{notation}[dummy]{Notation}

\newtheorem{trick}[dummy]{Calculation Device}
\begin{document}

\bibliographystyle{amsalpha}

\title{Graded Lie Algebras of Maximal Class II}

\date{17 February 1999}

\thanks{The first author is a member  of CNR-GNSAGA.  He
  acknowledges support from MURST-Italy, and from   CNR-GNSAGA
  under    \emph{Sottoprogetto ``Calcolo    Simbolico''  del  Progetto
    Strategico ``Applicazioni della Matematica  per la Tecnologia e la
    Societ\`a''.}  The first  author is grateful to SMS-ANU, Canberra,
  for  the kind hospitality  while  working at  this paper.  The first
  author is grateful to Giuseppe Jurman for help in the transition
  from  \AmS-\TeX\ to \amslatex/.\\
  The  second author  is grateful to   the Department of  Mathematics,
  University of Trento, for kind hospitality and support while working
  at this paper.}

\author{A.~Caranti}

\address{Dipartimento di Matematica\\
  Universit\`a degli Studi di Trento\\
  via Sommarive 14\\
  I-38050 Povo (Trento)\\
  Italy} 

\email{caranti@science.unitn.it} 

\urladdr{http://www-math.science.unitn.it/\~{ }caranti/}

\author{M.F.~Newman}

\address{School of Mathematical Sciences\\
  Australian National University\\
  Canberra, ACT 0200\\
  Australia }  

\email{newman@maths.anu.edu.au}  

\urladdr{http://wwwmaths.anu.edu.au/\~{ }newman/Mike\underline{ }Newman.html}



\begin{abstract}
  We describe the isomorphism classes  of infinite-dimensional  graded
  Lie algebras of maximal class over fields of odd characteristic.
\end{abstract}

\maketitle

\thispagestyle{empty}

\newcommand{\idgLams}[0]{infinite-dimensional graded Lie algebras of
maximal class}

\newcommand{\idgLam}[0]{infinite-dimensional graded Lie algebra of
maximal class}

\newcommand{\fdgLams}[0]{finite-dimensional graded Lie algebras of
maximal class}

\newcommand{\gLam}[0]{graded Lie algebra of maximal class}

\newcommand{\gLams}[0]{graded Lie algebras of maximal class}

\section{Introduction}

In this paper we describe a determination of
all isomorphism types of infinite-dimensional graded Lie algebras
with maximal class over fields having odd characteristic.

In our earlier paper with the same title~\cite{max}
we showed that there are
$|\F|^{\aleph_{0}}$ infinite-dimensional graded Lie algebras
of maximal class over a field $\F$ of positive characteristic
(Corollary 9.4).
This was done by describing the isomorphism type of such an algebra 
by its (normalized) sequence of two-step centralizers or, equivalently,
a sequence of elements from the projective line over $\F$.
We exhibited enough different sequences corresponding to algebras
by using a known construction,
which we called inflation, and combining it with taking limits.
The starting point was
some insoluble infinite-dimensional graded Lie algebras of maximal
class over the field with $p$ elements built by
A.~Shalev \cite{Sha:max} as positive parts of twisted loop algebras of
some finite-dimensional simple algebras of Albert and Frank with a
non-singular derivation on top. 

These AFS-algebras (see Section~\ref{sec:AFS}) together with
the process of inflation and taking limits
in fact suffice to allow us here to describe {\em all} \idgLams\ 
up to isomorphism. 
This required us to refine our techniques for determining which 
infinite sequences of elements of the projective line 
can arise from sequences of two-step centralizers of \idgLams. 

Carrara  has shown \cite{claretta, claretta2} that  each AFS-algebra $L$
can be characterised by a finite quotient, in the sense 
that the sequence of two-step centralizers of
AFS-algebras  are characterized by a certain initial  segment.
In other words, each AFS-algebra has a finite quotient which is not a
quotient of any other \gLam. We call the smallest such quotient $Q(L)$.
Our results show that for {\em odd} characteristic 
such  a statement fails for  all the other \idgLams.
This perhaps explains the very distinguished role that AFS-algebras play 
in the present context.

The above distinction is the basis for the proof of our main result
(Theorem~\ref{theorem:classification})  which states that
an \idgLam\ arises either via a finite number of inflation steps from
an AFS-algebra, or as the limit of an infinite number of inflation steps. 
(In the latter case
the resulting algebra is independent of the starting algebra.) 

A consequence of our  result is that every infinite-dimensional \gLam\ 
has    at    most   three    distinct    constituent   lengths    (see
Definition~\ref{def:constituents}).  This is not true in general though we
expect that an arbitrary \gLam\ has at most four constituent lengths.

Jurman~\cite{beppe2, beppe} has shown that the characteristic two case
is essentially (and not just technically) different by finding another
family of  uninflated \gLams.  Moreover  he has shown that,  with this
family       added,        a       result       corresponding       to
Theorem~\ref{theorem:classification}   holds.    His  proof   requires
further technical refinements.

There are finite-dimensional \gLams\ which do not arise as quotients
of \idgLams. These algebras are quotients of only finitely many other
\gLams.  Those which are not quotients of others are usually called
{\em terminal}. Carrara has determined all the terminal \gLams\ which
have a $Q(L)$ as a quotient for some AFS-algebra $L$. It might be
possible to refine her methods and our methods to obtain a description
of all terminal \gLams.

Computations with  the  {ANU} $p$-{Q}uotient  {P}rogram \cite{HNO} 
have been invaluable for understanding  the structure of the algebras 
under consideration, and  have  constantly   guided our  proofs.  
We   refer to~\cite{max} for more details. 
We note here just that we made use of
a more recent version of the program which allowed us occasionally
to study algebras with dimensions in the thousands to help clarify
matters.

We begin with some preliminaries: in Section~\ref{sec:background}
we review in particular the facts we need from~\cite{max}. In
Section~\ref{sec:leitfaden} we then state our
result, and give an outline of the proof. Although we use only
elementary methods, the proof is technically involved, so that we try
to state clearly the ideas behind the proof here. The actual proof is
carried out in Sections~\ref{sec:num_shorts}-\ref{sec:more_than_two}.

We refer to the previous paper as I.

We are grateful to Claretta Carrara and Giuseppe Jurman for reading
earlier versions of the manuscript and making useful suggestions.

\newpage

\section{Preliminaries}\label{sec:preliminaries}

Throughout the paper, except where explicitly stated otherwise, $p$ is
an odd prime, $\F$ is  a field of characteristic $p$, and
$\F_{p}$ denotes the field with $p$ elements.

We will use without further mention the Jacobi identity
\begin{equation*}
 [u [v w]] = [u v w] - [u w v];
\end{equation*}
unspecified brackettings are left-normed.
We will use many times the generalized
Jacobi identity
\begin{subequations}
  \begin{equation}\label{eq:Jacobi}
    [v [y \U{z}{\lambda}]] 
    =
    \sum_{i=0}^{\lambda} 
    (-1)^{i} \binom{\lambda}{i}
    [v \U{z}{i} y \U{z}{\lambda - i}],
  \end{equation}
  also in its equivalent form where we ``count from the back''
  \begin{equation}\label{eq:Jacobi_back}
    \begin{aligned}{}
      [v [y \U{z}{\lambda}]] 
      &=
      (-1)^{\lambda}
      \sum_{i=0}^{\lambda} 
      (-1)^{i} \binom{\lambda}{i}
      [v \U{z}{\lambda - i} y \U{z}{i}]
      \\&\uptoasign
      \sum_{i=0}^{\lambda} 
      (-1)^{i} \binom{\lambda}{i}
      [v \U{z}{\lambda - i} y \U{z}{i}],
    \end{aligned}
  \end{equation}
  where we have written  $a \uptoasign b$ to mean $a =  \pm b$.  As we
  will   be   applying~\eqref{eq:Jacobi_back}   only   when   $[v   [y
  \U{z}{\lambda}]]$ is  known to be zero,  the change of  sign will be
  irrelevant.
\end{subequations}

Lucas' Theorem (\cite{Lucas},  but see for instance \cite{Knuth:Lucas}
for some recent developments on  the subject), states that if $\lambda$
and $\mu$ are arbitrary non-negative integers, and
\begin{align*}
  \lambda &= \lambda_{0} + \lambda_{1} p + \dots + \lambda_{n} p^{n},\\
  \mu &= \mu_{0} + \mu_{1} p + \dots + \mu_{n} p^{n},
\end{align*}
with $0 \le \lambda_{i}, \mu_{i} < p$, are their \p-adic representations, then
\begin{equation*}
  \binom{\lambda}{\mu} \equiv 
  \prod_{i=0}^{n} \binom{\lambda_{i}}{\mu_{i}} \pmod{p}.
\end{equation*}
From a \p-adic point of view, if $q =  p^{h}$, and we write $\lambda =
\lambda_{0} + \lambda_{1}  q$  and $\mu = \mu_{0}  +  \mu_{1} q$, with
$\lambda_{0}, \mu_{0} < q$, then
\begin{equation*}
  \binom{\lambda}{\mu} \equiv 
  \binom{\lambda_{0}}{\mu_{0}} \cdot \binom{\lambda_{1}}{\mu_{1}} \pmod{p}.
\end{equation*}
We use    these   results many times   to   evaluate  binomial
coefficients modulo $p$.  For instance, if the remainder of $\lambda$ modulo
some power of $p$ is less than the corresponding remainder
of $\mu$, then $\dbinom{\lambda}{\mu} \equiv 0 \pmod{p}$.

We also use the following elementary facts: if $q$ is a power of $p$,
then
\begin{equation*}
  (-1)^{i} \dbinom{q-1}{i} \equiv 1 \pmod{p},
\end{equation*}
for $0  \le i \le q
- 1$;
and
\begin{equation}\label{eq:sum_of_binom}
 0 = (1 + (-1))^{n} = \sum_{i=0}^{n} (-1)^{i} \binom{n}{i}.
\end{equation}
A typical application of this is
$$
 \sum_{i=0}^{n-1} (-1)^{i} \binom{n}{i} = - (-1)^{n}.
$$

In several circumstances we will be using shorthands like
\begin{equation*}
  \lambda_{1}, (\lambda_{2}, \lambda_{3}^{\mu})^{\infty},
\end{equation*}
 where the $\lambda_{i}$ are arbitrary   objects and $\mu$ is a  non-negative
 integer, to denote the sequence
\begin{equation*}
 \lambda_{1}, 
 \lambda_{2}, \underbrace{\lambda_{3}, \dots, \lambda_{3}}_{\mu},
 \lambda_{2}, \underbrace{\lambda_{3}, \dots, \lambda_{3}}_{\mu},
 \dots \ .
\end{equation*}

We will introduce more specific notation and known facts in
Sections~\ref{sec:background} and \ref{sec:leitfaden}.

\section{Background}\label{sec:background}

This section recalls in  a concise way all the  material we   need
from~I, and introduces some notation  and conventions we use
throughout the paper.

\begin{definition}
An \emph{algebra of maximal class} is a graded Lie algebra 
\begin{equation*}
 L = \bigoplus_{i = 1}^{\infty} L_{i}
\end{equation*}
over a field $\F$ of characteristic $p$, such that $\dim (L_{1}) = 2$,
$\dim (L_{i}) \le 1$, for $i \ge  2$, and $L$ is generated by $L_{1}$. 
In other words, $[L_{i} L_{1}] = L_{i+1}$ for $i \ge 1$.
\end{definition}

\From now  on, let $L$ be  an infinite-dimensional  algebra of maximal
class over the field $\F$;  thus $\dim (L_{i}) =
1$ for all $i \ge 2$. The two-step centralizers
\begin{equation*}
 C_{i} = C_{L_{1}} (L_{i}),
\end{equation*}
for $i \ge  2$, are subspaces of dimension  $1$ of the two-dimensional
space $L_{1}$. We choose an element $y \in L_{1}$ so that
\begin{equation*}
C_{2} = \F y \ .
\end{equation*}
It follows from~I, Lemma~3.1, that  if $C_{i} =  \F y$ for all $i$,
then $L$  is isomorphic  to  the  infinite-dimensional
metabelian algebra of  maximal class $\Meta  (\F)$ over $\F$.  Suppose
then that  there is  a second \TSC\    and that $C_{n}$ is   its first
occurrence. We say that the second \TSC\ \emph{occurs  in weight $n$}. 
Write $C_{n} = \F x$ for some $x \in L_{1}$.  \From~I, Theorem~5.5, we
have
\begin{prop}\label{prop:parameter}
The second two-step centralizer first occurs in weight $2 p^{h}$,
for some $h \ge 1$.
\end{prop}

In  other words, $n = 2  p^{h}$, for some  $h$.  We will see in
Section~\ref{sec:leitfaden},
Step~\ref{step:occurrence_of_centralizers},   that a similar  statement
holds for all two-step centralizers.   If there are other centralizers
distinct  from $\F y$ and  $\F x$, then  $x$ and $y$  can be redefined so
that the third \TSC\ in order of occurrence is $\F \cdot (y - x)$. All
two-step centralizers can then be written uniquely as
\begin{equation*}
 C_{i} = \F \cdot (y - \alpha_{i} x),
\end{equation*}
for  some $\alpha_{i} \in  \F  \cup  \{  \infty \}$. (See~I,
Section~3.) Here $y - \infty \cdot x = x$.

In the rest of the paper we use without further mention the following 
\begin{notation}
  \label{def:xyz}
  Let  $L$ be  an   infinite-dimensional,  non-metabelian algebra   of
  maximal class. We  write $\F y = C_{2}$  for the first \TSC, and $\F
  x$ for the second one in order of occurrence. We also write $z = x +
  y$.
\end{notation}

The sequence  $( C_{i} )_{i \ge 2}$  determines  $L$ up to isomorphism
(I,  Theorem~3.2).  Equivalently,  $L$ is   determined by the
sequence  $( y -  \alpha_{i} x )_{i \ge  2}$.  We call either of these
two equivalent  sequences \emph{the sequence of two-step centralizers}
of $L$. We even mix usage. 
For  instance, $\Meta (\F)$ is  determined,  as recalled above, by the
constant sequence of two-step centralizers $(y, y, y, \dotsc )$.

Suppose  $L$ is   not isomorphic to  $\Meta  (  \F )$.   Because of~I,
Lemma~3.3,  the sequence of two-step  centralizers  of $L$ consists of
consecutive occurrences of $y$, interrupted by isolated occurrences of
different two-step centralizers.    Suppose that  in the sequence   of
two-step centralizers we have a pattern of the form
\begin{gather*}
 C_{i} \ne \F y\ , 
 \qquad
 C_{i+m} \ne \F y\ ,\\
 C_{i+1} = C_{i+2} = \dots = C_{i+m-1} = \F y\ .\\
\end{gather*}

\begin{definition}[Constituents]
  \label{def:constituents}
  We call the subsequence
  \begin{equation*}
    C_{i+1}, C_{i+2}, \dots, C_{i+m}
  \end{equation*}
  a  \emph{constituent}   of  $L$.    The  \emph{length}  of   such  a
  constituent  is the  number  $m$.   Note that  this  is a  different
  definition from that of~I, Section~3, where the final centralizer
  $C_{i+m}$ was  not included.  Experience  suggests this is  a better
  usage; we  apologise for the  inconvenience the change  causes.  It
  follows that the lengths of constituents in this paper are increased
  by one with respect to~I.
\end{definition}

Also, the following special definition turns out to be useful.
\begin{definition}[First constituent and parameter]
  \label{def:first_constituent}
  Suppose  the first occurrence of the  second two-step centralizer is
  $C_{2 q} = \F  x$. It is convenient to regard the subsequence
  \begin{equation*}
  C_{2}, C_{3}, \dots, C_{2 q}
  \end{equation*}
  as a  constituent (and thus the  first one) \emph{of length  $2 q$}. 
  We  call   $q$  the  \emph{parameter}   of  the  algebra   $L$.   By
  Proposition~\ref{prop:parameter}, we have $q = p^{h}$ for some $h$.
\end{definition}

It will be handy to have the following

\begin{definition}\label{def:end_of_constituent}
  A non-zero element $v$ of $L_{i+m}$  will be said  to be \emph{at the
    end} of the constituent
  \begin{equation*}
    C_{i+1}, C_{i+1}, \dots, C_{i+m} \ .
  \end{equation*}
  Similarly, a non-zero element $v$ of $L_{i}$ will be said to be
  \emph{at the beginning} of the constituent, or, equivalently,
  \emph{at the end} of the previous constituent.
\end{definition}

If $L$ has only the two distinct two-step  centralizers $\F y$ and $\F
x$, then the sequence of constituent lengths determines the algebra up
to isomorphism.

There are  strong  restrictions on  the  length  of  constituents, as
proved in~I, Proposition~5.6. In the  notation of the present paper,
this reads

\begin{prop}\label{prop:constituent_lengths}
  Let $L$ have parameter $q$. Then the constituent lengths take
  values in the set
  \begin{equation*}
    \{ 2 q \} \cup \{ 2 q - p^{\beta} : \text{for $0 \le \beta \le h$} \}.
  \end{equation*}
\end{prop}

In~I, Section~6,  a  procedure  called  \emph{inflation} is
defined, which given an  algebra $L$ of maximal  class, and a subspace
$M = \F   w$ of  $L_{1}$ of  dimension  $1$, yields   another  algebra
$\Infl{M}{L}$ of maximal class. There is a  way of identifying $L_{1}$
and  $(\Infl{M}{L})_{1}$ (see~I,  Section~6).   Once this  is
done,   one  can describe  the  sequence  of two-step  centralizers of
$\Infl{M}{L}$,    starting   with    $C_{(\Infl{M}{L})_{1}}      \big(
(\Infl{M}{L})_{2} \big)$, as
\begin{equation*}
  \underbrace{w, \dots, w}_{2 p - 2}, 
  C_{2}, 
  \underbrace{w, \dots, w}_{p - 1}, 
  C_{3}, 
  \dots, 
  C_{i}, 
  \underbrace{w, \dots, w}_{p - 1}, 
  C_{i+1}, 
  \dotsc
  \ .
\end{equation*}
In   other words, inflation  can be  described by  inserting $p-1$
occurrences of  $w$ between every  pair of  consecutive terms of  the
sequence of two-step centralizers of  $L$ to obtain the corresponding
sequence for $\Infl{M}{L}$. The initial part of the sequence requires
special treatment.

There  is another procedure,  called \emph{deflation}  (I, Section~7),
that is a one-sided inverse of inflation. It associates to $L$ another
algebra of maximal class $L^{\downarrow}$.   Here, too, there is a way
of  identifying   $L_{1}$  with  $(L^{\downarrow})_{1}$.    With  this
identification, the sequence of two-step centralizers $(C_{i}')_{i \ge
  2}$  of the  deflation $L^{\downarrow}$  of $L$  is  described in~I,
Proposition~7.1, which we recall here as

\begin{prop}\label{prop:deflation}
  If for $i \ge 2$ all the two-step centralizers
  \begin{equation*}
  C_{i p}, C_{i p + 1}, \dots, C_{(i+1) p - 1}
  \end{equation*}
  coincide  with $C_{2} =  \F y$, then  $C_{i}' = \F  y$. If one of
  them equals $\F u \ne \F y$, then $C_{i}' = \F u$.
\end{prop}

It       should       be       noted       that       because       of
Proposition~\ref{prop:constituent_lengths},  at  most  one  among  the
listed two-step centralizers can differ from $\F y$.

It is clear from the above that we always have
\begin{equation*}
(\Infl{M}{L})^{\downarrow} = L.
\end{equation*}
In  other words  if an  algebra of  maximal class  is an  inflation of
another  algebra, then  the latter  algebra  is the  deflation of  the
former.

\begin{definition}
  An algebra $L$ of maximal class is said to be \emph{inflated} if
  \begin{equation*}
    L = \Infl{M}{(L^{\downarrow})}
  \end{equation*}
  for a suitable maximal subspace $M$ of $(L^{\downarrow})_{1}$.
\end{definition}

The Albert-Frank-Shalev algebras and their limits (Section~4
of~I) are not inflated.
This follows from the fact that inflated algebras can be explicitly
characterized (I, Proposition~7.4). In the notation of the current paper, this
can be formulated as

\begin{prop}\label{prop:inflation}
  An algebra $L$ is inflated if and only if all its constituents have
  length a multiple of $p$. In this case we have
  \begin{equation*}
    \Infl{M}{(L^{\downarrow})} = L,
  \end{equation*}
  where $M = C_{L_{1}} (L_{2})$ in the identification of $L_{1}$ with
  $(L^{\downarrow})_{1}$\ .
\end{prop}

There is a particular case which is easy to describe.
\begin{prop}\label{prop:inflation_2}
  Inflation  with  respect  to  the  first \TSC\  has  the  effect  of
  multiplying all  constituent lengths  by $p$. The  resulting algebra
  has  parameter greater than  $p$, and  all constituents  lengths are
  multiples of $p$.
  
  Conversely, suppose $L$ has parameter  $q$ greater than $p$, and all
  constituents have  length a multiple  of $p$.  Then the  sequence of
  constituent lengths of $L^{\downarrow}$ is obtained from that of $L$
  by dividing all constituent lengths by $p$.
\end{prop}

There is a special case of repeated  deflation that is particularly
useful. It follows immediately from Proposition~\ref{prop:deflation}.

\begin{prop}\label{prop:deflation_by_q}
  Suppose $L$  has parameter $q = p^{h}$.  The result of deflating $h$
  times  is an  algebra  in  which $x$  plays   the role  of the first
  centralizer. If all two-step centralizers
  \begin{equation*}
    C_{i q}, C_{i q + 1}, \dots, C_{(i+1) q - 1}
  \end{equation*}
  coincide with $C_{2} =  \F y$, then  $C_{i}' = \F  y$ in the algebra
  deflated $h$ times.  If  one of them equals  $\F  u \ne \F  y$, then
  $C_{i}' = \F u$ in the algebra deflated $h$ times.
\end{prop}

Here, too, Proposition~\ref{prop:constituent_lengths} implies that at most one
among the listed two-step centralizers can differ from $\F y$.

The following Lemma is proved in \cite{CaJur}. We will give an independent
proof in Section~\ref{sec:more_than_two}.

\begin{lemma}\label{lemma:third_after_short}
  A  two-step   centralizer other  than   the first  two  in  order of
  occurrence always occurs at  the end of  a short constituent, and it
  is followed by another short constituent.
\end{lemma}

Once the characterization is obtained, one can read off from it a much
more precise statement.

\section{Leitfaden}\label{sec:leitfaden}

In this section we state our main result,  and describe the pattern of
the proof.

Let $p$ be  a fixed odd prime.  Let  $\F$ be a field of
characteristic  $p$.    Consider   the   following two    classes   of
infinite-dimensional algebras of maximal class, defined over the field
$\F_{p}$ with $p$ elements:

\begin{itemize}
\item $\AFS (a, b, n, p)$, the insoluble Albert-Frank-Shalev algebras,
  for  $0 < a  <  b \le n$,   as  described in  \cite{Sha:max}, and~I,
  Section~4.
\item   $\AFS (a,  b, \infty, p)$,  the  soluble   limits  associated to  the
  Albert-Frank-Shalev algebras,  for $0 < a  < b$, as  described in~I,
  Section~9.
\end{itemize}

Note that all these algebras have two distinct two-step centralizers. 

In~I  we made  a different  (equivalent) choice  of parameters  in the
$\AFS$-algebras  as $0  \le a  <  b< n$;  the current  choice has  the
advantage  of  offering  a  more uniform  description.  (Compare  also
Proposition~7.3 of~I.) What we called $\AFS  (0, b, n, p)$ in~I is now
called $\AFS (b, n, n, p)$.  In the notation of the present paper, the
algebra $\AFS (a, b, n, p)$ has sequence of constituent lengths
\begin{equation*}
  2 q, (q^{r-2}, 2 q - 1, (q^{r-2}, 2 q)^{s - 1})^{\infty},
\end{equation*}
where $q =  p^{a}$, $r = p^{b-a}$,  and $s = p^{n-b}$. As mentioned
in Section~\ref{sec:preliminaries},  the notation $(q^{r-2}, 2 q)^{s -
  1}$ denotes $s -  1$ repetitions of  the pattern $q^{r-2}, 2 q$, and
$\pi^{\infty}$, where  $\pi$  is a   pattern of  constituent  lengths,
denotes  periodic repetition of  $\pi$.  Here we  understand of course
that when $n  = b$, so that $s  - 1 = p^{n-b}  -  1 = 0$,  the pattern
$(q^{r-2}, 2 q)$ does not occur at all.

In the notation of the present paper, the algebra $\AFS(a,b,\infty,p)$
has \SOCL{}
\begin{equation*}
 2 q, q^{r-2}, 2 q - 1, (q^{r-2}, 2 q)^{\infty}.
\end{equation*}
This explains our notation.

We have the following

\begin{theorem}[Characterization]\label{theorem:classification}
  Let $L$ be an infinite-dimensional algebra of maximal class over the
  field $\F$ of odd characteristic $p$.

  Then $L$ is obtained
  \begin{itemize}
  \item  either via a  finite number (possibly  zero)  of inflation steps
    from one of the algebras $\AFS (a,  b, n, p) \otimes_{\F_{p}} \F$, 
  \item or as the limit of an infinite number of inflation steps.
  \end{itemize}
\end{theorem}

We note that  when one inflates  an algebra infinitely many times, the
result  is independent of the algebra  one began with,  so that is not
necessary to specify the latter in this case.

In  the rest of  the paper, we take $L$  to be an infinite-dimensional
algebra of maximal class over a fixed field $\F$ of odd characteristic
$p$.

The key  to our  proof is that  we show that  \emph{every non-inflated
  algebra  of maximal class  is an  $\AFS$-algebra.}  Suppose  we have
done this, and let $L$ be any algebra of maximal class.  If $L$ is not
inflated,  we are finished.   If $L$  is inflated,  we note  its first
two-step   centralizer  $\kappa_{1}$,   and  consider   its  deflation
$L^{\downarrow}$.  If  $L^{\downarrow}$ is not inflated, it  has to be
an      $\AFS$-algebra,      and      we      recover      $L$      as
$\Infl{\kappa_{1}}{L^{\downarrow}}$.  If $L^{\downarrow}$ is inflated,
we note its first  two-step centralizer $\kappa_{2}$, and consider its
deflation  $(L^{\downarrow})^{\downarrow}$.   We  keep repeating  this
procedure. If after a finite number $m$ of deflation steps we reach an
algebra $N$ that is not  inflated, then $N$ will be an $\AFS$-algebra,
and  we  will  recover $L$  from  $N$  by  inflating with  respect  to
$\kappa_{m}, \dots,  \kappa_{2}, \kappa_{1}$.  If all  the algebras we
obtain in  the process are inflated,  then $L$ can be  obtained as the
limit of  an infinite  number of inflation  steps with respect  to the
sequence $(\kappa_{i})_{i \ge 1}$.

Our proof of Theorem~\ref{theorem:classification} is given in a number
of  steps. We have  first,  from  Proposition~\ref{prop:parameter} and
Proposition~\ref{prop:constituent_lengths},

\begin{step}\label{step:constituent_lengths}
  If $L$  is not metabelian, then  the first constituent has length $2
  q$,  for some power  $q = p^{h}$  of $p$. All other constituents
  have length $2 q$ or $2 q - p^{\beta}$, for $0 \le \beta \le h$.
\end{step}

\From now on we will assume $L$  is not metabelian,  and we will write
$q = p^{h}$ for its parameter.

We  call a constituent of length  $2 q$ {\em  long}, one of length $q$
{\em  short}.  All  others are {\em  intermediate}.

\From Proposition~\ref{prop:inflation} we have

\begin{step}\label{step:inflated}
  $L$ is  inflated if   and  only if all  constituents  have  length a
  multiple of $p$.
\end{step}

Step~\ref{step:inflated} clearly has the special case

\begin{step}\label{step:all_short}
  If all constituents  after the  first one  are  short,  then $L$  is
  inflated.
\end{step}

Because  of Step~\ref{step:all_short},  we  may assume  that there  is
another non-short constituent after the first long one. We first deal,
roughly  speaking, with  the case  when  there are  only two  distinct
two-step centralizers.

We begin with proving

\begin{step}\label{step:Longs_and_Shorts}
  Suppose  $L$ has two distinct  two-step  centralizers. If the second
  non-short constituent  is long, then  all constituents are  short or
  long, and thus $L$ is inflated.
\end{step}

This  follows  from  the  results  of Section~\ref{sec:num_shorts}~and
\ref{sec:Longs_and_Shorts}.

We have to consider the case when the second non-short constituent is
intermediate.

\begin{step}\label{step:shorts_before_intermediate}
  Suppose the second non-short  constituent is an intermediate one, of
  length  $2 q  -  p^{\beta}$, for  $0 \le  \beta  < h$,  and that  no
  two-step centralizer  higher than the  second one occurs before  it. 
  Then the initial pattern of constituents is of the form
  \begin{equation*}
    2 q, q^{r - 2}, 2 q - p^{\beta},
  \end{equation*}
  for some power $r = p^{k}$ of $p$.
\end{step}

This is proved in Section~\ref{sec:num_shorts}.

Note that such a constituent  pattern appears when one inflates
an $\AFS$-algebra  $\beta$ times with respect to the first
two-step centralizer.

\From now on we use the notation $r = p^{k}$.

Most  of the  next step  has already  been obtained  by  C.~Carrara in
\cite{claretta,     claretta2}.    We     include    a     proof    in
Section~\ref{sec:Step_6} for completeness.

\begin{step}\label{step:pattern_for_AFS}
  Suppose there is an initial segment of the \SOTC\ of $L$ that
  involves only the first two \TSCs, and that the corresponding
  segment of the \SOCL\ is
  \begin{equation*}
    2 q, q^{r - 2}, 2 q - p^{\beta},
  \end{equation*}
  where $0 \le \beta < h$. Then the algebra  has two distinct two-step
  centralizers, and the \SOCL\ consists of repetitions of the patterns
  \begin{equation*}
    2 q, q^{r - 2} \qquad \text{or} \qquad  2 q - p^{\beta}, q^{r - 2}.
  \end{equation*}
\end{step}

The following is clear in view of Proposition~\ref{prop:inflation_2}.

\begin{step}\label{step:AFSl}
  Suppose $L$ has two distinct two-step centralizers.  If the sequence
  of constituent lengths is
  \begin{equation*}
    2 q, q^{r-2}, 2 p - p^{\beta}, (q^{r-2}, 2 q)^{\infty},
  \end{equation*}
  then  $L$ is obtained  from an  $\AFS(a,b,\infty,p) \otimes_{\F_{p}}
  \F$ via $\beta$  inflation steps with respect to  the first two-step
  centralizer.
\end{step}

We thus suppose that there is another intermediate constituent. In
Section~\ref{sec:AFS} we prove

\begin{step}\label{step:getting_to_Claretta}
  Suppose there is an initial segment of the \SOTC\ of $L$ that
  involves only the first two \TSCs, and that the corresponding
  segment of the \SOCL\ is
  \begin{equation*}
    2 q, q^{r-2}, 2 q - p^{\beta}, (q^{r-2}, 2 q)^{s-1}, 
         q^{r-2}, 2 q - p^{\beta}.
  \end{equation*}
  Then $s = p^{l}$ for some ${l}$.
\end{step}

Again, this is consistent with the structure of repeated inflations of
the $\AFS$-algebras with respect to the first two-step centralizer.

\From now on we will use the notation $s = p^{{l}}$.

We now appeal to the main result of \cite{claretta}.

\begin{theorem}[C.~Carrara]\label{claretta}
  Let $L$ be an infinite-dimensional algebra of maximal class over the
  field $\F_{p}$, where $p$ is any prime.
  
  Suppose the \SOTC\ has an initial segment in which only two distinct
  two-step centralizers occur,  and the corresponding  initial segment
  of the \SOCL\ is
  \begin{equation*}
    2 q, q^{r-2}, 2 q - 1, (q^{r-2}, 2 q)^{s-1}, q^{r-2}, 2 q - 1,
  \end{equation*}
  with $q = p^{h}$, $r = p^{k}$ and $s = p^{l}$.
  
  Then $L \cong \AFS (h, h + k, h + k + l, p)$.  In particular, $L$ has
  two distinct \TSCs.
\end{theorem}

Using this, we immediately obtain

\begin{step}
  Suppose the \SOTC\ of $L$ has  an initial segment involving only the
  first two two-step  centralizers, and that the  corresponding \SOCL\ 
  begins with
\begin{equation*}
 2 q, q^{r-2}, 2 q - p^{\beta}, (q^{r-2}, 2 q)^{s-1}, q^{r-2}, 2 q - p^{\beta}.
\end{equation*}
Then  $L$ has two distinct  two-step centralizers, and $L$ is obtained
from an $\AFS$-algebra via $\beta$ inflation steps with respect to the
first two-step centralizer.
\end{step}

We  now deal  with  the case  of more  than  two distinct  \TSCs.  The
following        critical    result               is      proved    in
Section~\ref{sec:specialization}.

\begin{prop}[Specialization of two-step centralizers]
\label{prop:specialization}
  Let $L$ be an infinite-dimensional algebra of maximal class with
  three or more distinct two-step centralizers.
  
  Let $y$ be the first  two-step centralizer and $w$ a two-step
  centralizer distinct from $y$.

  Choose elements $x_{1}$ and $x_{2}$ in one of the following ways

  \begin{center}
    $x_{1} = y$ and $x_{2} = w$
    \quad
    or 
    \quad
    $x_{1} = w$ and $x_{2} = y$.
  \end{center}
  
  Consider the sequence of two-step centralizers of $L$. Leave $x_{1}$
  fixed, and change  all two-step centralizers other than  $y$ and $w$
  to $x_{2}$.
  
  Then the resulting sequence is the sequence of two-step centralizers
  of    an  algebra  of maximal   class   with  two  distinct two-step
  centralizers.
\end{prop}

We   use    Proposition~\ref{prop:specialization}     to        extend
Proposition~\ref{prop:parameter}.

\begin{step}\label{step:occurrence_of_centralizers}
  The $i$-th two-step centralizer  in order of occurrence appears  for
  the first time as
  \begin{equation*}
    C_{L_{1}} (L_{2 p^{n}}),
  \end{equation*}
  for some $n \ge i - 1$.
\end{step}

This is proved in Section~\ref{sec:more_than_two}.

Finally we obtain

\begin{step}\label{step:more_than_two_is_inflated}
  Suppose $L$ has at least three distinct two-step centralizers.  Then
  all     constituents  are either  short      or  long,  so  that  by
  Step~\emph{\ref{step:inflated}} the algebra is inflated.
\end{step}

This is proved in Section~\ref{sec:more_than_two}, and completes the proof.

\section{Sequences of short constituents}\label{sec:num_shorts}

In~I we have proved that the second constituent is short,
provided the characteristic of the field is odd.  Here we will
investigate more generally how many short constituents one gets
between two non-short ones in odd characteristic.  In doing this, we
will introduce a few basic arguments that we use repeatedly in
the course of the paper.  These will be presented in some
detail the first time they are used.

Because of Step~\ref{step:all_short}, we may assume that not all
constituents after the first one are short.

If the second non-short constituent is long, and no \TSC\ other than
the first two occurs before it, then deflation quickly tells us
how many short constituents occur between the two long ones. In fact,
suppose the \SOCL\ of $L$ begins as
\begin{equation*}
  2 q, q^{m}, 2 q.
\end{equation*}
The sequence of two-step centralizers therefore begins, starting
with $C_{2 q}$, as
\begin{equation*}
x, {\underbrace{y, \dots, y}_{q-1}, x, \underbrace{y, \dots,
    y}_{q-1}, x, \dots, \underbrace{y, \dots, y}_{q-1}, x},
\underbrace{y, \dots, y}_{2 q - 1},
\end{equation*}
where the pattern $\underbrace{y, \dots, y}_{q-1}, x$ appears $m$
times.  If we deflate $h$ times (remember that we use the standing
notation $q = p^{h}$), according to
Proposition~\ref{prop:deflation_by_q}, we obtain that the deflated
algebra has sequence of two-step centralizers, beginning with
$C'_{2}$,
\begin{equation*}
  \underbrace{x, \dots, x}_{m+1}, y \ .
\end{equation*}
We obtain  that the  length of the  first constituent in  the deflated
algebra   is,  according   to  Definition~\ref{def:first_constituent},
$m+3$.   By Step~\ref{step:constituent_lengths}, $m  + 3  = 2  r$, for
some power $r = p^{k}$ of $p$, so that $m = 2 r - 3$. We have obtained

\begin{lemma}\label{lemma:shorts_before_long}
  Suppose the \SOTC\ of $L$ has an initial segment involving only the
  first two two-step centralizers, and that the corresponding \SOCL\ 
  is
  \begin{equation*}
    2 q, q^{m}, 2 q.
  \end{equation*}
  Then $m = 2 r - 3$ for some power $r = p^{k}$ of $p$.
\end{lemma}

If the  second non-short constituent is intermediate,  the situation is
more  complicated.   Consider  for  instance  the algebra  $L  =  \AFS
(1,2,2,3)$, which  has two distinct two-step  centralizers, and \SOCL\ 
$6, (3,  5)^{\infty}$, so that the  parameter $q$ is $3$  here. We now
write down  the initial segment of  the \SOTC\ of  $L$, beginning with
$C_{2 q} = C_{6}$, isolating the two-step centralizers in groups of $q
= 3$. We get
\begin{equation*}
  x y y | x y y | y y x | y y x | y y y | y x y | y x y | y y y | 
  x y y | \dots \ .
\end{equation*}
The resulting \SOTC\ in $L^{\downarrow}$ is, according to
Proposition~\ref{prop:deflation},
\begin{equation*}
  x, x, x, x, y, x, x, y, x, \dots \ .
\end{equation*}
Here the intermediate constituent 
does not provide the end of the first constituent
in the  deflated algebra. We will  need considerably more  work in the
rest of this  section to deal with this case. 

It  might be worth noting
that  if  all  constituents  are  long or  short,  then  all  two-step
centralizers other than the first one occur in the form $C_{i}$, where
$i \equiv 0 \pmod{q}$.   The occurrence of an intermediate constituent
disturbs   this  pattern,   and  it   takes  further   occurrences  of
intermediate  constituents to  bring it  back into  step.  However, it
takes  an intermediate  constituent of  maximal length  $2 q  -  1$ to
disrupt the fact  that all two-step centralizers other  than the first
one occur in the form $C_{i}$, where $i \equiv 0 \pmod{p}$.

We now introduce a simple but important

\begin{trick}
  Let $0 \ne v \in L_{i}$, for some $i > 1$. Suppose we know that all
  the two-step centralizers $C_{i}, C_{i+1}, \dots, C_{i+n-1}$ equal
  one of the first two, that is, they are either $y$ or $x$.
  Consider, according to Definition~\ref{def:xyz}, the element $z = x
  + y \in L_{1}$, a notation that we will use in the rest of the
  paper.  We have that if $[v x_{1} x_{2} \dots x_{n}]$ is a non-zero
  commutator in $L$, with $x_{i} \in \{ x, y \}$, then
  \begin{equation*}
    [v x_{1} x_{2} \dots x_{n}] = [v \U{z}{n}] \ .
  \end{equation*}
\end{trick}

We begin with a simple generalization of Lemma~5.3 of~I.

\begin{lemma}\label{fact:odd}
  Suppose no \TSC\  other than the first two  occurs before the second
  non-short  constituent.   Then  the  number  of  short  constituents
  between the first constituent and the second non-short one is odd.
\end{lemma}

\begin{proof}
  Suppose the sequence of constituent lengths starts with
  $
    2 q, q^{m},
  $
  where  $m$ is even.
  We will prove that this is followed by another
  short constituent.
   
  We have to prove
  \begin{equation}\label{eq:odd_num_shorts}
    [y \U{z}{q (m + 3) - 1} x] = 0.
  \end{equation}
  Since $q (m + 3) - 1$ is
  even, we can  consider the integer
  \begin{equation*}
    \lambda = \dfrac{q (m + 3) -1}{2} = q \dfrac{m + 2}{2} +
    \dfrac{q-1}{2} \ .
  \end{equation*}
  We will obtain equation~\eqref{eq:odd_num_shorts} from the
  expansion of the following expression with the generalized Jacobi
  identity~\eqref{eq:Jacobi_back}:
  \begin{equation} \label{eq:expansion}
    \begin{aligned}{}
      0 &= [ [ y \U{z}{\lambda} ] [ y \U{z}{\lambda} ] ] \\&\uptoasign
      \sum_{i=0}^{\lambda} (-1)^{i} \binom{\lambda}{i} [y \U{z}{2 \lambda
        - i} y \U{z}{i}] .
    \end{aligned}
  \end{equation}
  Note that the total weight of the commutators in~\eqref{eq:expansion}
  is that of the left-hand side of~\eqref{eq:odd_num_shorts}:
  \begin{equation*}
    2 \lambda + 2 = q (m + 2) + q - 1 + 2 = q (m + 3) + 1.
  \end{equation*}
  
  We now  look at the  resulting terms in  \eqref{eq:expansion}. Since
  $y$ is the first centralizer, the term $[ y \U{z}{2 \lambda - i} y]$
  will vanish  most of the  time, except when  $[ y \U{z}{2  \lambda -
    i}]$  occurs in  a weight  that marks  the end  of a  constituent. 
  Keeping  in mind  the information  we already  have on  the previous
  constituents, the sum reduces to
\begin{equation}\label{eq:summ_for_odd_num_of_shorts}
 \begin{aligned}{}
   0 &= [y \U{z}{2 \lambda} y] + \\&\phantom{=\ } + \Bigg( \sum_{j =
     1}^{(m + 2)/2} (-1)^{q j} \binom{q \frac{m + 2}{2} +
     \frac{q-1}{2}}{q j} \Bigg) \cdot [y \U{z}{2 \lambda + 1}] .
 \end{aligned}
\end{equation}
Lucas' Theorem yields
\begin{equation*}
  (-1)^{q j} \binom{q \frac{m + 2}{2} + \frac{q-1}{2}}{q j} =
  (-1)^{j} \binom{\frac{m + 2}{2}}{j};
\end{equation*}
by~\eqref{eq:sum_of_binom}
\begin{equation*}
\sum_{j = 1}^{(m+2)/2} (-1)^{j} \binom{\frac{m+2}{2}}{j} = - 1.
\end{equation*}
Therefore \eqref{eq:summ_for_odd_num_of_shorts} yields
\begin{equation*}
0 = [y \U{z}{2 \lambda} z] - [y \U{z}{2 \lambda} y] = [y \U{z}{2
  \lambda} x],
\end{equation*}
as $z = x + y$. This is precisely \eqref{eq:odd_num_shorts}.
\end{proof}

\begin{lemma}\label{fact:shorts_after_intermediate_1}
  Suppose $L$ has an initial segment of the \SOTC\ involving only the
  first two \TSCs, and suppose the corresponding segment of the \SOCL\ 
  is
  \begin{equation*}2 q, q^{m}.\end{equation*}
  
  Suppose further that an intermediate constituent in the \SOTC\ of
  $L$ is \emph{preceded} by at least $m$ short ones ending in $x$.
  
  Then that intermediate constituent is
  \emph{followed} by at least $m$ short ones ending in $x$.
\end{lemma}

We will make a corresponding statement later
(Lemma~\ref{fact:shorts_after_long}).

\begin{proof}
  Let the intermediate constituent have length $2 q - p^{\beta}$.
  Suppose we have already proved that this intermediate constituent is
  followed by $\tau$ short ones ending in $x$, with $1 \le \tau < m$,
  so that there is a segment
  \begin{equation*}
  \dotsc, q^{m}, 2 q - p^{\beta}, q^{\tau}
  \end{equation*}
  in the \SOCL\ of $L$.  We will show that the next constituent is
  also short and ends in $x$. Let $\sigma = m - 1 - \tau \ge 0$.
  
  According      to  Definition~\ref{def:end_of_constituent}, take   a
  homogeneous   element $v$  at   the   beginning of the   $\sigma$-th
  constituent before the intermediate one.  What has to be proved is
  \begin{equation}\label{eq:3_stgh}
    [v \U{z}{q (m+2) - p^{\beta}} x] = 0.
  \end{equation}
  We assume this is not the case, and obtain a contradiction.

  We expand
  \begin{equation*}
    0 = [v [y \U{z}{q (m+2) - 1} x]] \ ;
  \end{equation*}
  here $[y \U{z}{q (m+2)  - 1} x] =  0$ because of our assumption on
  the initial segment  of the \SOTC. Note  that this will give us terms
  of   weight     $p^{\beta}$    higher   than     the   left-hand term
  of~\eqref{eq:3_stgh}.
  
We get, using~\eqref{eq:Jacobi} and the fact that $[v x] = 0$,
\begin{align*}
 0
 &=
 [v [y \U{z}{q (m+2) - 1} x]]
 \\&=
 [v [y \U{z}{q (m+2) - 1}] x]
 \\&=
  [v \U{z}{q (m+2)} x] \cdot
 \Bigg(
  \sum_{i = 0}^{\sigma} (-1)^{q i} \binom{q (m + 1) + q - 1}{q i}
 \\&\qquad
  +
  \sum_{i = \sigma + 2}^{m+1} (-1)^{q (i - 1) + q - p^{\beta}} 
      \binom{q (m + 1) + q - 1}{q (i - 1) + q - p^{\beta}}
 \Bigg)
 \\&\phantom{=\ }
 +
  (-1)^{q (m+2) - p^{\beta}} \binom{q (m+2) - 1}{q (m+2) - p^{\beta}}
  [v \U{z}{q (m+2) - p^{\beta}} y \U{z}{p^{\beta}-1} x]
 \\&=
 [v \U{z}{q (m+2)} x] \cdot \sum_{i=0}^{m} (-1)^{i} \binom{m+1}{i}
 \\&\qquad+ 
 (-1)^{m+1} [v \U{z}{q (m+2) - p^{\beta}} y \U{z}{p^{\beta}-1} x]
 \\&=
 (-1)^{m+1} \cdot
 \big(- [v \U{z}{q (m+2) - p^{\beta}} z \U{z}{p^{\beta}-1} x] 
 \\&\qquad+
 [v \U{z}{q (m+2) - p^{\beta}} y \U{z}{p^{\beta}-1} x]\big)
 \\&\uptoasign
 [v \U{z}{q (m+2) - p^{\beta}} x \U{z}{p^{\beta}-1} x] .
\end{align*}
Here we have used several times Lucas' Theorem and its consequences.

Now if~\eqref{eq:3_stgh} does not  hold we obtain a  contradiction. In
fact~\eqref{eq:3_stgh} can  fail in  two  ways. It might be   that $[v
\U{z}{q (m+2) -  p^{\beta}}]$ is centralized  by $y$. This  means that
the $(\tau + 1)$-th  constituent  after the  intermediate one  is  not
short, so by Proposition~\ref{prop:constituent_lengths} it has length
at least
\begin{equation*}
2 q - \frac{q}{p} = q + (p-1) \frac{q}{p} > q + \frac{q}{p},
\end{equation*}
as $p > 2$. However, we would have just obtained a constituent of length
$q + p^{\beta} \le q + q/p$, a contradiction.

As an alternative, it could be  that $[v \U{z}{q (m+2) - p^{\beta}} w]
= 0$ for some  $w$, with $\F w \ne \F y, \F x$.   We would then have a
short  constituent not  ending in  $x$; but  we have  shown  that this
constituent would be followed by  a constituent of length $p^{\beta} <
q$, and this is a final contradiction.
\end{proof}

\begin{lemma}\label{lemma:initial_num_of_shorts}
  Suppose  no \TSC\ other than the  first two occurs before the second
  non-short constituent.  Suppose the \SOCL\  of $L$  begins as $2  q,
  q^{m}, 2 q - p^{\beta}$.
  
  Then $m = p^{k} - 2$ for some $k$.
\end{lemma}

This yields in particular Step~\ref{step:shorts_before_intermediate}.

\begin{proof}
  Deflating         $\beta$         times,        according         to
  Proposition~\ref{prop:inflation_2}, we  may assume $\beta  = 0$.
  
  Write $m + 2  = r n$, where $r = p^{k}$ is  a power of $p$ (possibly
  $1$),      and     $n      \not\equiv      0     \pmod{p}$.       By
  Lemma~\ref{fact:shorts_after_intermediate_1}, we know that the first
  intermediate  constituent is  followed by  at least  $m$  short ones
  ending in $x$.  Let
  \begin{equation*}
    \mu = q \big( m + 4 + r - 1 \big) = q \big( r (n + 1) + 1 \big) .
  \end{equation*}
  If $n > 1$, we have $m \ge r - 1$.  Therefore we obtain
  \begin{equation*}
    0 \ne [y \U{z}{\mu - 2} y] = [y \U{z}{\mu - 2} z] =
    [y \U{z}{\mu - 1}] .
  \end{equation*}
  We will prove
  \begin{equation}\label{eq:to_be_proved}
    0 = [y \U{z}{\mu}] = [y \U{z}{\mu - 2} y z],
  \end{equation}
  so that  the  algebra  is  finite-dimensional,  since
  constituents have length greater than $1$. This contradiction will give
  that $n = 1$, so that $m = p^{k} - 2$ as claimed.

  We consider
  \begin{equation*}
    \text{$\lambda = \eta + q \theta$, where $\eta = \dfrac{q-1}{2}$, and
      $\theta = r \dfrac{n+1}{2}$.}
  \end{equation*}
  Here $n$ is odd, because of Lemma~\ref{fact:odd}. Note that
  \begin{equation*}
    2 \lambda + 2 = \mu + 1,
  \end{equation*}
  so that $[ [y \U{z}{\lambda}] [y  \U{z}{\lambda}] ]$ has the same
  weight as the commutators in~\eqref{eq:to_be_proved}.

  We expand, using~\eqref{eq:Jacobi_back} and the information about the
  constituents,
  \begin{equation}
    \label{eq:two_summations}
    \begin{aligned}
      0
      &=
      [ [y \U{z}{\lambda}] [y \U{z}{\lambda}] ]
      \\&\uptoasign
      [y \U{z}{\mu}] \cdot
      \Bigg(
      \sum_{i = 0}^{r - 1} (-1)^{1 + q i} \binom{\eta + q \theta}{1 + q i}
      \\&\qquad+
      \sum_{i = r + 1}^{\theta} (-1)^{q i} \binom{\eta + q \theta}{q i}
      \Bigg) \ .
    \end{aligned}
  \end{equation}
  
  Now
  \begin{align*}
    \sum_{i = 0}^{r-1} (-1)^{1 + q i} \binom{\eta + q \theta}{1 + q i}=
    &=
    - \eta \sum_{i = 0}^{r-1} (-1)^{i} \binom{\theta}{i}
    \\&=
    - \eta \binom{\theta}{0} 
    \\&= 
    - \eta
    \\&\equiv
    \frac{1}{2} \pmod{p},
  \end{align*}
  as $\dbinom{r (n+1)/2}{i} \equiv 0 \pmod{p}$ for $0 < i < r$.
  As to the second summation, we have
  \begin{align*}
    \sum_{i = r + 1}^{\theta} (-1)^{q i} \binom{\eta + q \theta}{q i}
    &=
    \sum_{i = r + 1}^{\theta} (-1)^{i} \binom{\theta}{i}
    \\&=
    \sum_{j = 2}^{(n+1)/2} (-1)^{j} \binom{\frac{n+1}{2}}{j}
    \\&=
    - \left( 1 - \frac{n+1}{2} \right).
  \end{align*}
  The total coefficient  of $[y \U{z}{\mu}]$ is thus
  \begin{equation*}
    \frac{1}{2} - 1 + \frac{n+1}{2} = \frac{n}{2} \not\equiv 0 \pmod{p}, 
  \end{equation*}
  as claimed.
\end{proof}
 
\begin{lemma}\label{fact:shorts_after_long}
  Suppose the \SOTC\ of $L$ has  an initial segment involving only the
  first two two-step  centralizers, and that  the corresponding \SOCL\ 
  is
  \begin{equation*}2 q,  q^{r  - 2},\end{equation*}
  where  $r$  is a power  of $p$.  Then every long
  constituent  in the \SOTC\  of $L$ is followed by  at  least $r - 2$
  short ones ending in $x$.
\end{lemma}

\begin{proof}
  By  Lemma~\ref{fact:shorts_after_intermediate_1}   and induction, we
  may  assume   that we have  the  subsequence  $\dotsc, q^{r-2}, 2 q,
  q^{\tau}$ in the \SOTC\ of $L$, for some $0 \le \tau < r - 2$. Let
  $\sigma = r - 2 - \tau$, so that $\sigma \le r  - 2$, and let $v$ be
  a  homogeneous element  at  the end   of  the  $\sigma$-th  two-step
  centralizer before the long one, so that $[v x] = 0$. We have to prove
  \begin{equation*}
  [v \U{z}{q r} x] = 0.
  \end{equation*}
  
  We  have $[ y \U{z}{q r  -  1} x] = 0$,   from our assumption on the
  initial segment of the \SOTC\ of $L$. We compute with~\eqref{eq:Jacobi}
  \begin{align*}
    0
    &=
    [v [y \U{z}{q r - 1} x] ]
    \\&=
    [v [y \U{z}{q r - 1}] x ]
    \\&=
    [v \U{z}{q r} x] \cdot \Bigg(
    \sum_{i = 0}^{\sigma-1} (-1)^{q i} \binom{q r - 1}{q i} +
    \sum_{i = \sigma+1}^{r-1} (-1)^{q i} \binom{q r - 1}{q i}
    \Bigg)
    \\&=
    [v \U{z}{q r} x] \cdot \Bigg(
    - (-1)^{\sigma} \binom{r-1}{\sigma}
    \Bigg)
    \\&=
    - [v \U{z}{q r} x]. &&\qed
  \end{align*}
  \renewcommand{\qed}{}
\end{proof}

\section{Algebras with long and short constituents}
\label{sec:Longs_and_Shorts}

In this section we prove Step~\ref{step:Longs_and_Shorts}. We consider
algebras of maximal class  $L$ with two distinct two-step centralizers
in which  the second non-short constituent  is long, and show that all
constituents             are        short     and           long.   By
Lemma~\ref{lemma:shorts_before_long}, we   know that the   sequence of
constituent lengths begins like
$$
 2 q, q^{2 r - 3}, 2 q,
$$
where $q = p^{h}$ and $r = p^{k}$. The algebra we obtain by deflating
$h$ times has thus first constituent of length $2 r$, as in the proof
of Lemma~\ref{lemma:shorts_before_long}.
 
By Lemma~\ref{fact:shorts_after_long} and
induction, we assume that at some point in the sequence of constituent
lengths we have the pattern 
\begin{equation}\label{eq:long_plus_shorts}
 \dotsc, q^{r-2}, 2 q, q^{r-2}.
\end{equation}
We first show that no  intermediate constituent can now occur. Suppose
there is  a constituent of length $2  q - p^{\beta}$, for  some $0 \le
\beta  <   h$,  immediately  following   \eqref{eq:long_plus_shorts}.  
Deflating $\beta$ times, we may assume $\beta$ to be zero.

We use the relation
$$
 [y \U{z}{q (r + 1) -1} x] = 0,
$$
commuted against a homogeneous element $v \ne 0$ of suitable weight, so that
the resulting commutator 
\begin{equation}
  \label{eq:to_be_killed}
  [v \U{z}{q (r + 1)} y]
\end{equation}
falls at the end of the final constituent in
$$
 \dots, q^{r-2}, 2 q, q^{r-2}, 2 q - 1.
$$
Inspection shows that $v$ falls within the long constituent in
\eqref{eq:long_plus_shorts}; in particular, $[v y] = 0$.

We have,
using \eqref{eq:Jacobi_back}, 
\begin{align*}
 0
 &=
 [v [y \U{z}{q (r + 1) -1} x] ]
 \\&=
 [v [y \U{z}{q (r + 1) -1}] x]
 -
 [v x [y \U{z}{q (r + 1) -1}] ]
 \\&\uptoasign
 [v \U{z}{q (r + 1)} x]
 \cdot
 \Bigg(
 \sum_{i=0}^{r-2} (-1)^{2 q - 2 + q i} 
  \binom{q (r + 1) -1}{2 q - 2 + q i}
 \Bigg)
 \\&\phantom{=\ }
 - [v \U{z}{q (r + 1)} y]
 \\&\phantom{=\ }
 - [v \U{z}{q (r + 1)} z]
  \Bigg(
 \sum_{i=0}^{r-2} (-1)^{2 q - 1 + q i} 
  \binom{q (r + 1) -1}{2 q - 1 + q i}
 \Bigg).
\end{align*}
Now
\begin{align*}
 \sum_{i=0}^{r-2} (-1)^{2 q - 2 + q i} 
  \binom{q (r + 1) -1}{2 q - 2 + q i}
 &=
 \sum_{i=0}^{r-2} (-1)^{2 q - 2 + q i} 
  \binom{q r + q - 1}{q (i + 1) + q - 2}
 \\&=
 \sum_{i=0}^{r-2} (-1)^{i} 
  \binom{r}{i + 1}
 \\&=
 - \sum_{i=1}^{r-1} (-1)^{i} 
  \binom{r}{i}
 \\&= - (1 + (-1)^{r}) = 0 \ ,
\end{align*}
and similarly
\begin{equation*}
 \sum_{i=0}^{r-2} (-1)^{2 q - 1 + q i} 
  \binom{q (r + 1) -1}{2 q - 1 + q i}
 =
 0 \ .
\end{equation*}
We obtain that~\eqref{eq:to_be_killed}  vanishes,  thereby providing a
contradiction.

We now show that if there is a further short constituent at the end of 
\eqref{eq:long_plus_shorts}, then this is followed by another short one.
Using the appropriate $v$, which is again centralized by $y$, we compute with
\eqref{eq:Jacobi_back}  
\begin{align*}
 0
 &=
 [ v [y \U{z}{q (r + 1) - 1} x] ]
 \\&=
 [v [y \U{z}{q (r + 1) -1}] x]
 -
 [v x [y \U{z}{q (r + 1) -1}] ]
 \\&=
  [v \U{z}{q (r + 1)} x]
 \cdot
 \Bigg(
 \sum_{i=0}^{r-1} (-1)^{q - 1 + q i} 
  \binom{q (r + 1) -1}{q - 1 + q i}
 \Bigg)
 \\&\phantom{=\ }
 - [v \U{z}{q (r + 1)} y]
 \\&\phantom{=\ }
 - [v \U{z}{q (r + 1)} z]
  \Bigg(
 \sum_{i=1}^{r} (-1)^{q i} 
  \binom{q (r + 1) -1}{q i}
  \Bigg).
\end{align*}
We have easily
\begin{align*}
 \sum_{i=0}^{r-1} (-1)^{q - 1 + q i} 
  \binom{q (r + 1) -1}{q - 1 + q i} 
 &= 
  \sum_{i=0}^{r-1} (-1)^{i} 
  \binom{r}{i} 
 \\&=
 - (-1)^{r} =
 1,
\intertext{and similarly}
 \sum_{i=1}^{r} (-1)^{q i} 
  \binom{q (r + 1) -1}{q i}
 &=
 -1,
\end{align*}
so that
\begin{equation*}
  0
  =
   [v \U{z}{q (r + 1)} x]
  - [v \U{z}{q (r + 1)} y]
  + [v \U{z}{q (r + 1)} z]
 = 2 [v \U{z}{q (r + 1)} x],
\end{equation*}
as required.

We finally show that no pattern
\begin{equation}\label{eq:too_many}
 \dots, q^{r-2}, 2 q, q^{\sigma}, 2 q - 1
\end{equation}
is allowed for $\sigma > r  - 1$. (Because  of the deflation argument
we  have  already  employed   at the   beginning  of   the   proof  of
Lemma~\ref{lemma:initial_num_of_shorts},  this  takes   care of    any
trailing  intermediate  constituent here.)   In  fact, this  should be
followed by  Lemma~\ref{fact:shorts_after_intermediate_1} by at least
$q^{r-2}$.  Deflating  $h$ times, we see by  direct  inspection of the
sequence of   two-step centralizers,   very   much as we     did after
Lemma~\ref{lemma:shorts_before_long},  that we obtain a constituent of
length at least
\begin{equation*}
 \sigma + 2 + r - 2 + 1 > 2 r.
\end{equation*}
We have noted at the beginning of this section that the deflated algebra has
parameter $2 r$; we have thus obtained a contradiction.

\section{Proof of Step~\ref{step:pattern_for_AFS}}\label{sec:Step_6}

We are assuming that the \SOTC\ of the algebra $L$ has an initial segment
involving only two distinct two-step centralizers, and that the corresponding
\SOCL\ is
\begin{equation}\label{eq:AFS_branch}
 2 q, q^{r-2}, 2 q - p^{\beta},
\end{equation}
for some $0 \le \beta < h$.
We want to prove that the whole algebra $L$ has two distinct two-step
centralizers, and that the sequence of constituent lengths consists of
repetitions of the patterns
\begin{equation*}
 2 q, q^{r-2} 
 \qquad 
 \text{or} 
 \qquad
 2 q - p^{\beta}, q^{r-2} \ .
\end{equation*}

As mentioned in the Introduction, most of the  proofs of this section,
except for minor technical differences, and for one point we will
comment upon later, were obtained first by C.~Carrara in her PhD
Thesis~\cite{claretta}. 

From         Lemmas~\ref{fact:shorts_after_intermediate_1}         and
\ref{fact:shorts_after_long},    and   induction,    we    know   that
under~\eqref{eq:AFS_branch}   there   are   at   least   $r-2$   short
constituents  ending in  $x$ after  a non-short  one. So  we  start by
assuming, by induction, that we have a segment
\begin{equation}\label{eq:a}
 \dots, q^{r-2}, 2 q, q^{r-2},
\end{equation}
where all short constituents end in  $x$, and show first that the next
constituent  is  \emph{not}  short.  This implies  in  particular,  by
Lemma~\ref{lemma:third_after_short},   that    no   further   two-step
centralizer can occur. We use the  relation $[y \U{z}{(r+1)q - 1} y] =
0$, which  says that in~\eqref{eq:AFS_branch} the  last constituent is
\emph{not} short. Let  $v$ be an element at the end  of the last short
constituent before the  long one in~\eqref{eq:a}. We have  to show $[v
\U{z}{(r+1)q} y] = 0$. We compute with~\eqref{eq:Jacobi_back}
\begin{align*}
 0
 &=
 [v [y \U{z}{(r+1)q - 1} y]]
 \\&=
 [v [y \U{z}{(r+1)q - 1}] y]
 -
 [v y [y \U{z}{(r+1)q - 1}] ] \ .
\end{align*}
The first term expands, up to a sign, as
\begin{multline*}
 [v \U{z}{(r+1)q} y] \cdot
 \Bigg(
 \sum_{i=0}^{r-2} (-1)^{q-1 + q i} \binom{q - 1 + q r}{q-1 + q i}
 +
 (-1)^{q-1 + q r} \binom{q - 1 + q r}{q - 1 + q r}
 \Bigg)
 \\ -
 [v \U{z}{(r+1)q} y] ,
\end{multline*}
while the second one expands, up to the same sign, to
\begin{equation*}
 - [v \U{z}{(r+1)q} z] \cdot
 \Bigg(
 \sum_{i=1}^{r-1} (-1)^{i} \binom{q - 1 + q r}{q i}
 \Bigg) .
\end{equation*}
Summing up, we obtain
\begin{align*}
 0
 &=
 [v \U{z}{(r+1)q} y] \cdot
 \big( - (-1)^{r} + (-1)^{r} - 1 \big)
 \\&\phantom{=\ }
 + [v \U{z}{(r+1)q} z] \cdot \big( 1 + (-1)^{r} \big)
 \\&=
 - [v \U{z}{(r+1)q} y].
\end{align*}

An entirely similar argument shows  that no constituent of the form $2
q - p^{\gamma}$, with $\gamma \not\in \{ \beta, h \}$ can occur at the
end  of~\eqref{eq:a}.  We use  the  same  $v$,  and the  relation  $[y
\U{z}{(r+2) q  - p^{\gamma} - 1} y]  = 0$, which states  that the last
constituent in~\eqref{eq:AFS_branch}  is \emph{not}  of length $2  q -
p^{\gamma}$. We have,  using~\eqref{eq:Jacobi_back} as in the previous
case,
\begin{align*}
 0
 &=
 [v [ y \U{z}{(r+2) q - p^{\gamma} - 1} y]]
 \\&=
 [v [ y \U{z}{(r+2) q - p^{\gamma} - 1}] y]
 -
 [v y [ y \U{z}{(r+2) q - p^{\gamma} - 1}] ]
 \\&\uptoasign
 [v  \U{z}{(r+2) q - p^{\gamma}} y] \cdot
 \\&\qquad\cdot
 \Bigg(
 \sum_{i=1}^{r-1} (-1)^{q i + q - p^{\gamma} - 1} 
 \binom{q (r+1) + q - p^{\gamma} - 1}{q i + q - p^{\gamma} - 1}
 \\&\qquad\phantom{\cdot \Bigg(\ }
 + (-1)^{q (r+1) + q - p^{\gamma} - 1}
  \Bigg)
 \\&\phantom{=\ }
 - [v  \U{z}{(r+2) q - p^{\gamma}} y]
 \\&\phantom{=\ }
 - 
 [v  \U{z}{(r+2) q - p^{\gamma}} z] \cdot
 \\&\qquad\cdot
 \Bigg( 
 \sum_{i=1}^{r} (-1)^{q i + q - p^{\gamma}} 
 \binom{ q (r+1) + q - p^{\gamma} - 1}{q i + q - p^{\gamma}}
 \Bigg)
 \\&=
 \big( - (-1)^{q - p^{\gamma} - 1}(1 + (-1)^{r}) - 1 \big) 
 [v  \U{z}{(r+2) q - p^{\gamma}} y]
 \\&=
 - [v  \U{z}{(r+2) q - p^{\gamma}} y].
\end{align*}
Here we  have  used the  fact that $\gamma  \ne h$,  so that  $0 < q -
p^{\gamma} - 1  < q -  p^{\gamma} < q$,  and $\dbinom{q - p^{\gamma} -
  1}{q - p^{\gamma}} = 0$.

We now assume, again proceeding
by induction, that we have a segment
\begin{equation}\label{eq:b}
 \dots, q^{r-2}, 2 q - p^{\beta}, q^{r-2},
\end{equation}
and show that the next constituent can be either long or of length $2 q -
p^{\beta}$. In both cases, by Lemma~\ref{lemma:third_after_short}, no further two-step
centralizer can thus occur. We distinguish three cases.

Suppose first, by way of contradiction, that in \eqref{eq:b} a short constituent
follows. We can deflate $\beta$ times, according to Proposition~\ref{prop:inflation_2}, and
suppose that
\eqref{eq:b} takes the form
\begin{equation}\label{eq:b_not_another_short}
 \dots, q^{r-2}, 2 q - 1, q^{r-1}.
\end{equation}
Take $v$  to be an  element at the  end of the last  short constituent
before  the  intermediate one.  In  particular,  $[v  y] \ne  0$.  The
configuration~\eqref{eq:b_not_another_short}   implies   that   $[   v
\U{z}{\lambda} y] \ne 0$, where
\begin{equation*}
  \lambda = q (r + 1) - 1 = q r + q - 1. 
\end{equation*}
We will prove $[ v \U{z}{\lambda} y z] = [ v \U{z}{\lambda + 2}] =
0$, a contradiction, as  we also have $[ v \U{z}{\lambda} y  y] = [ v
\U{z}{\lambda + 1} y] = 0$, since constituents have positive length.

We use  the  relation $[y  \U{z}{\lambda} y] =   0$, stating that  the
\SOCL\  does \emph{not}   start   as  $2   q,    q^{r-1}$.   Expanding
with~\eqref{eq:Jacobi_back}, we get,    using   the fact  that   $[  v
\U{z}{\lambda} y y] = [ v \U{z}{\lambda + 1} y] = 0$, as above,
\begin{align*}
 0
 &=
 [v [y \U{z}{\lambda} y] ]
 \\&\uptoasign
 [v y [y \U{z}{\lambda}] ] 
 \\&\uptoasign
 [v \U{z}{\lambda + 2}] \cdot
  \Bigg(
   \sum_{i = 0}^{r-1} (-1)^{q i + 1} \binom{q r + q - 1}{q i + 1}
  \Bigg)
 \\&=
 [v \U{z}{\lambda + 2}].
\end{align*}

Now we deal with the case when~\eqref{eq:b} continues as
\begin{equation}\label{eq:different_intermediate}
 \dots, q^{r-2}, 2 q - p^{\beta}, q^{r-2}, 2 q - p^{\gamma},
\end{equation}
with $h > \gamma > \beta$. Deflating $\beta$ times, according to
Proposition~\ref{prop:inflation_2}, we may assume
$\beta$ to be zero. We employ a variation of the previous argument. We use the same $v$, so
we assume
$[ v \U{z}{\lambda} y] \ne 0$,
where this time
\begin{equation*}
 \lambda = q (r + 1) + q - p^{\gamma} - 1, 
\end{equation*}
and prove $[ v \U{z}{\lambda + 2}] = 0$, a contradiction. We use the relation $[y
\U{z}{\lambda} y] = 0$,  stating that the \SOCL\ does \emph{not} start as $2 q, q^{r-2},
2 q - p^{\gamma}$. 
Expanding with~\eqref{eq:Jacobi_back}, we get
\begin{align*}
 0
 &=
 [v [y \U{z}{\lambda} y] ]
 \\&\uptoasign
 [v y [y \U{z}{\lambda}] ]
 \\&\uptoasign
 [v \U{z}{\lambda + 2}] \cdot
 \\&\qquad\cdot
  \Bigg(
   (-1)^{1} \binom{q (r + 1) + q - p^{\gamma} - 1}{1}
   \\&\qquad\phantom{\cdot \Bigg( \ }+
   \sum_{i = 1}^{r-1} (-1)^{q i + q - p^{\gamma}} 
   \binom{q (r + 1) + q - p^{\gamma} - 1}{q i + q - p^{\gamma}}
  \Bigg)
 \\&=
 [v \U{z}{\lambda + 2}],
\end{align*}
as
\begin{equation*}
 \binom{q (r + 1) + q - p^{\gamma} - 1}{q i + q - p^{\gamma}}
 \equiv
 \binom{q (r + 1)}{q i} \cdot \binom{q - p^{\gamma} - 1}{q - p^{\gamma}},
\end{equation*}
since $0 < q - p^{\gamma} - 1 < q - p^{\gamma} < q$.

Finally, we consider the case   when~\eqref{eq:different_intermediate}
holds   for  some  $\gamma  <    \beta$.  This   has  no   counterpart
in~\cite{claretta}, where one works under the  assumption $\beta = 0$. 
This    time  we   may    deflate    $\gamma$  times,  according    to
Proposition~\ref{prop:inflation_2}, and assume $\gamma$ to be zero, so
that~\eqref{eq:b} takes the form
\begin{equation}\label{eq:ba}
 \dots, q^{r-2}, 2 q - p^{\beta}, q^{r-2}, 2 q - 1.
\end{equation}
We start
with the same $v$ as before, that is, at the end of the last short constituent before the
one of length $2 q - p^{\beta}$. The pattern \eqref{eq:ba} amounts to say
$[v \U{z}{\lambda} y] \ne 0$, where
\begin{equation*}
 \lambda = q (r + 1) + q - p^{\beta} - 1.
\end{equation*}
We will prove 
\begin{equation*}
 0
 = 
 [v \U{z}{\lambda + p^{\beta} + 1}] 
 = 
 [v \U{z}{\lambda} y \U{z}{p^{\beta}}].
\end{equation*}
This would exhibit a constituent of length at most $p^{\beta} < q$, a
contradiction.

We define $u = [v \U{z}{p^{\beta} - 1}]$, so that in particular $[u y] = 0$. We have to prove
\begin{equation*}
 [u \U{z}{\lambda + 2}] = 0.
\end{equation*}
We use the relation
$[y \U{z}{\lambda} x] = 0$ which says that the second non-short constituent is of length $2
q - p^{\beta}$. We expand from the back, according to~\eqref{eq:Jacobi_back},
\begin{align*}
 0
 &=
 [ u [y \U{z}{\lambda} x] ]
 \\&=
 [ u [y \U{z}{\lambda}] x ]
 -
 [ u x [y \U{z}{\lambda}]]
 \\&\uptoasign
 [u \U{z}{\lambda + 1} x] \cdot
 \\&\qquad\cdot
 \Bigg(
  (-1)^{p^{\beta} - 1} \binom{\lambda}{p^{\beta}-1}
  +
  \sum_{i=2}^{r} (-1)^{q i + p^{\beta} - 2} 
   \binom{\lambda}{q i + p^{\beta} - 2}
 \Bigg)
 \\&\phantom{=\ }
 - [u \U{z}{\lambda + 1} z] \cdot
 \Bigg(
 (-1)^{p^{\beta}} \binom{\lambda}{p^{\beta}}
 +
  \sum_{i=2}^{r} (-1)^{q i + p^{\beta} - 1} 
   \binom{\lambda}{q i + p^{\beta} - 1}
 \Bigg)
 \\&\phantom{=\ }
 - [u \U{z}{\lambda + 1} y].
\end{align*}

We evaluate the binomial coefficients modulo $p$, according to Lucas' Theorem.
\begin{align*}
 (-1)^{p^{\beta} - 1} \binom{\lambda}{p^{\beta}-1}
 &=
 (-1)^{p^{\beta} - 1} \binom{q (r + 1) + q - p^{\beta} - 1}{p^{\beta}-1}
 \\&\equiv
 (-1)^{p^{\beta} - 1} \binom{q -  2 p^{\beta} + p^{\beta} - 1}{p^{\beta}-1}
 \\&\equiv
 \binom{p^{\beta} (p^{h - \beta} -  2)}{0}
 (-1)^{p^{\beta} - 1} \binom{p^{\beta} - 1}{p^{\beta}-1}
 \\&\equiv
 1.
\intertext{In a similar fashion}
 (-1)^{q i + p^{\beta} - 2} 
   \binom{\lambda}{q i + p^{\beta} - 2}
  &\equiv
  (-1)^{i} \binom{r + 1}{i} \cdot (-1)^{p^{\beta}-2}
  \binom{q -  2 p^{\beta} + p^{\beta} - 1}{p^{\beta}-2}
  \\&\equiv
  (-1)^{i} \binom{r + 1}{i},
\end{align*}
and
\begin{equation*}
  (-1)^{q i + p^{\beta} - 1} 
   \binom{\lambda}{q i + p^{\beta} - 1}
   \equiv
   (-1)^{i} \binom{r + 1}{i}.
\end{equation*}
Also
\begin{align*}
   (-1)^{p^{\beta}} \binom{\lambda}{p^{\beta}}
   &\equiv
   (-1)^{p^{\beta}}
   \binom{q -  2 p^{\beta} + p^{\beta}}{p^{\beta}}
   \\&\equiv
   (-1)^{p^{\beta}} 
   \binom{p^{\beta} (p^{h-\beta} - 2)}{p^{\beta}}
   \\&\equiv
   2.
\end{align*}
Keeping in mind that
\begin{equation*}
 \sum_{i = 2}^{r} (-1)^{i} \binom{r + 1}{i}
 =
 - \Bigg( 1 - \binom{r+1}{1} + \binom{r+1}{r+1} \Bigg)
 =
 - 1,
\end{equation*}
and that $[u \U{z}{\lambda+1} y] = 0$,
we obtain
\begin{align*}
 0
 &=
 [ u [y \U{z}{\lambda} x] ]
 \\&=
 [u \U{z}{\lambda+1} x] \cdot 
 \Bigg( 
  1 + \sum_{i = 2}^{r} (-1)^{i} \binom{r + 1}{i}
 \Bigg)
 \\&\phantom{=\ }
 - [u \U{z}{\lambda+1} z] \cdot
 \Bigg(
  2 + \sum_{i = 2}^{r} (-1)^{i} \binom{r + 1}{i}
 \Bigg)
 \\&\phantom{=\ }
 \\&\uptoasign
 [u \U{z}{\lambda+1} z],
\end{align*}
as required.

\section{The algebras of Albert-Frank-Shalev}\label{sec:AFS}

We now prove Step~\ref{step:getting_to_Claretta}.

The situation  is the  following.  We have  an  algebra $L$
with initial segment of the \SOCL\ of the form
\begin{equation*}
 2 q, q^{r - 2}, 2 q - p^{\beta},
\end{equation*}
for some $0 \le \beta < h$, where no two-step centralizer other than
the first two occurs. Here $q = p^{h}$ and $r = p^{k}$ as usual.
We know  from Step~\ref{step:pattern_for_AFS} that the \SOCL\ consists
of repetitions of the patterns
\begin{equation*}
  2 q, q^{r - 2} \qquad \text{or} \qquad  2 q - p^{\beta}, q^{r - 2}.
\end{equation*}
In view  of Proposition~\ref{prop:inflation_2}, we may deflate $\beta$
times,  and assume that  $\beta = 0$. We can  also suppose, in view of
Step~\ref{step:AFSl}, that another intermediate constituent appears in
the \SOCL. The \SOCL\ of $L$ thus begins like
\begin{equation}\label{eq:pattern_for_AFS.1}
 2 q, q^{r-2}, 2 q - 1, (q^{r-2}, 2 q)^{\sigma - 1}, q^{r - 2}, 2 q - 1.
\end{equation}
We have to show that $\sigma = p^{l}$ for some $l \ge 0$.

We  first show  that   $\sigma$ is odd.   \eqref{eq:pattern_for_AFS.1}
implies that  the commutator $[y \U{z}{\mu} y]$,  of weight $\mu + 2$,
where $\mu = q r ( \sigma + 1)  + 2 q -  3$, is non-zero. We will prove
that  if    $\sigma$ is   even,  then   $[y  \U{z}{\mu} y]    =  0$, a
contradiction. Let
\begin{equation*}
 \lambda = q r \frac{\sigma}{2} + q \frac{r + 1}{2} + \frac{q - 3}{2},
\end{equation*}
so that $2 \lambda = \mu$. We expand with~\eqref{eq:Jacobi_back}
\begin{align*}
 0
 &=
 [ [ y \U{z}{\lambda} ] [ y \U{z}{\lambda} ] ]
 \\&\uptoasign
 [y \U{z}{\mu} y] + \theta [y \U{z}{\mu + 1}].
\end{align*}
Now the first term of the coefficient $\theta$ is, up to a sign,
\begin{equation*}
 \binom{\lambda}{2 q - 1}
 =
 \binom{q r \frac{\sigma}{2} + q \frac{r + 1}{2} + \frac{q - 3}{2}}{q + q -1}
 \equiv 0 \pmod{p},
\end{equation*}
by Lucas'  Theorem, as  the top entry  in the binomial  coefficient is
congruent  to $\dfrac{q  - 3}{2}  \pmod{q}$, while  the bottom  one is
congruent to $q  - 1 > \dfrac{q  - 3}{2}$. It is not  difficult to see
that all terms in $\theta$  vanish similarly.  Therefore $\theta = 0$,
and $[y \U{z}{\mu} y] = 0$, as claimed.

We now prove that the \SOCL\ of $L$ cannot have an initial
segment of the form
\begin{equation}
  \label{eq:pattern_for_AFS.2}
  \begin{aligned}
    2 q, &q^{r-2}, 2 q - 1, (q^{r-2}, 2 q)^{\sigma - 1},\\
         &q^{r - 2}, 2 q - 1,  (q^{r-2}, 2 q)^{\tau - 1},\\
         &q^{r-2}, 2 q - 1, q,
  \end{aligned}
\end{equation}
for   any  $\tau   <  \sigma$.   If  there   is  a   segment   of  the
form~\eqref{eq:pattern_for_AFS.2}, then  the commutator $[y \U{z}{\mu}
y]$, of weight $\mu + 2$, where
\begin{equation*}
 \mu = q r (\sigma + \tau + 1) + 3 q - 4,
\end{equation*}
is  non-zero. We  will  show that  this  commutator vanishes,  thereby
obtaining a contradiction. We use the relation $[y \U{z}{\lambda} y] =
0$, where
\begin{equation*}
 \lambda = q r ( \tau + 1 ) + 2 q - 3.
\end{equation*}
This  relation  states that  the  initial  segment  of the  \SOCL\  is
\emph{not} of the form
\begin{equation*}
 2 q, q^{r-2}, 2 q - 1, (q^{r-2}, 2 q)^{\tau - 1}, q^{r - 2}, 2 q - 1.
\end{equation*}
(Remember  we have  taken $\tau <  \sigma$.) Take  an  appropriate
non-zero homogeneous element  $v$ so  that $[v [y \U{z}{\lambda}  y]]$ has
the appropriate weight $\mu +  2$. This means $v$  has weight $\mu  -
\lambda = q r \sigma + q - 1$. Inspection shows that $v$ lies in the
middle   of the    last    long  constituent  in   the   first   line
of~\eqref{eq:pattern_for_AFS.2}. Therefore $[v \U{z}{q}]$  lies at the end of
this long  constituent, and of  course $[v  y]  =  0$.  As  we  start
expanding with~\eqref{eq:Jacobi}
\begin{equation*}
  0
 =
 [v [y \U{z}{\lambda} y]]
 =
 [v [y \U{z}{\lambda}] y],
\end{equation*}
we first encounter a coefficient
\begin{equation*}
 \sum_{i = 1}^{r - 1} (-1)^{q i} \binom{\lambda}{q i},
\end{equation*}
as we  go through the  $r - 2$  short constituents following  the long
one. This evaluates modulo $p$, according to Lucas' Theorem, to
\begin{align*}
 \sum_{i = 1}^{r - 1} (-1)^{q i} \binom{q r ( \tau + 1 ) + q + q - 3}{q i}
 &\equiv
 \sum_{i = 1}^{r - 1} (-1)^{i} \binom{r ( \tau + 1 ) + 1}{i}
 \\&\equiv
 (-1)^{1} \binom{r ( \tau + 1 ) + 1}{1}
 \\&\equiv
 - 1.
\end{align*}
We  claim that     all   the     remaining  binomial      coefficients
in~\eqref{eq:Jacobi}  vanish. In fact the  next one,  as we go through
the   first    intermediate     constituent   in    the    second line
of~\eqref{eq:pattern_for_AFS.2}, is
\begin{equation*}
 \binom{q r ( \tau + 1 ) + q + q - 3}{q r + q - 1}.
\end{equation*}
This vanishes, because the top entry is congruent to $q - 3 \pmod{p}$,
while the bottom one is congruent to $q - 1 > q - 3$. The same applies
to all  the remaining coefficients;  as we pass the  last intermediate
constituent,  the  bottom entry in  the  binomial coefficient  becomes
congruent to $q - 2 > q - 3$.

Let now  $p^{l}$ be the highest  power of  $p$ (possibly $1$) dividing
$\sigma$. Write $\sigma = p^{l}  \sigma'$, so that $\sigma' \not\equiv
0 \pmod{p}$. We suppose $p^{l} \ne \sigma$, so  that $p^{l} < \sigma$,
and obtain a contradiction.

In fact we have shown that  \eqref{eq:pattern_for_AFS.2} does not hold
for $\tau \le p^{l}$. By Step~\ref{step:pattern_for_AFS}, the
\SOCL\ of $L$ must have an initial segment of the form
\begin{equation}\label{eq:pattern_for_AFS.3}
 2 q, q^{r-2}, 2 q - 1, (q^{r-2}, 2 q)^{\sigma - 1}, q^{r - 2}, 2 q - 1,
 (q^{r-2}, 2 q)^{p^{l}}.
\end{equation}
In particular, the commutator $[y \U{z}{\mu} y]$, of weight $\mu + 2$, where
\begin{equation*}
 \mu = q r (\sigma + p^{l} + 1) + 2 q - 3,
\end{equation*}
is  non-zero. We  will  show that  this  commutator vanishes,  thereby
obtaining a contradiction. Take
\begin{equation*}
 \lambda = q r \frac{\sigma + p^{l}}{2} + q \frac{r + 1}{2} + \frac{q - 3}{2},
\end{equation*}
so that $2 \lambda = \mu$, and
expand with~\eqref{eq:Jacobi_back} the balanced relation
\begin{equation}
  \label{eq:summation}
  \begin{aligned}
    0
    &=
    [ [y \U{z}{\lambda}] [y \U{z}{\lambda}] ]
    \\&\uptoasign
    [y \U{z}{\mu} y]
    \\&
    \phantom{=\ }
    + [y \U{z}{\mu + 1}] \cdot
    \Bigg(
    \sum_{i = 0}^{p^{l} - 1} \sum_{j = 2}^{r} 
    (-1)^{q r i + q j} \binom{\lambda}{q r i + q j}
    \Bigg) .
  \end{aligned}
\end{equation}
A word  of comment is in order.  The summation in~\eqref{eq:summation}
comprises  the terms  of~\eqref{eq:Jacobi_back} as  we go  through the
final         segment         $(q^{r-2},         2         q)^{p^{l}}$
of~\eqref{eq:pattern_for_AFS.3}.  As we  pass the  second intermediate
constituent  in~\eqref{eq:pattern_for_AFS.3}, the  bottom entry  of the
relevant binomial  coefficient becomes congruent to $q  - 1 \pmod{q}$,
while the top  entry $\lambda$ is congruent to $\dfrac{q -  3}{2} < q -
1$. Therefore the terms we have not displayed vanish as above.

We compute the summation in~\eqref{eq:summation}. Note first
\begin{equation*}
  \sum_{i = 0}^{p^{l} - 1} \sum_{j = 2}^{r} 
  (-1)^{q r i + q j} \binom{\lambda}{q r i + q j}
  =
  \sum_{i = 0}^{p^{l} - 1} (-1)^{i}
  \sum_{j = 2}^{r} 
  (-1)^{j} \binom{\lambda}{q r i + q j}.
\end{equation*}
For a given $i$ we now have, using Lucas' Theorem,
\begin{equation}
  \label{eq:i}
  \begin{aligned}
    (-1)^{i} \sum_{j = 2}^{r} 
    &(-1)^{j} \binom{\lambda}{q r i + q j} =
    \\&=
    (-1)^{i} \sum_{j = 2}^{r} 
    (-1)^{j} 
    \begin{pmatrix}
      &q r \frac{\sigma + p^{l}}{2} &+ &q \frac{r + 1}{2} &+ &\frac{q - 3}{2}\\
      &q r i &+ &q j & &
    \end{pmatrix}
    \\&=
    (-1)^{i}   \binom{\frac{\sigma + p^{l}}{2}}{i}
    \sum_{j = 2}^{r-1} 
    (-1)^{j} 
    \binom{\frac{r + 1}{2}}{j}
    +
    (-1)^{i+1} \binom{\frac{\sigma + p^{l}}{2}}{i+1}.
  \end{aligned}
\end{equation}
Now $(r+1)/2 \le r-1$, as $r \ge 3$. Therefore
\begin{align*}
  \sum_{j = 2}^{r-1} 
  (-1)^{j} 
  \binom{\frac{r + 1}{2}}{j}
  &=
  \sum_{j = 2}^{(r+1)/2} 
  (-1)^{j} 
  \binom{\frac{r + 1}{2}}{j} 
  \\&=
  - \left( 1 - \frac{r + 1}{2} \right) 
  \\&\equiv
  - \frac{1}{2} \pmod{p}.
\end{align*}
Suppose first $p^{l} \ne 1$. As $p^{l}$ divides $\sigma$, and $i <
p^{l}$, Lucas' Theorem yields that 
$\dbinom{\frac{\sigma + p^{l}}{2}}{i} \not\equiv 0 \pmod{p}$ only for
$i = 0$, and $\dbinom{\frac{\sigma + p^{l}}{2}}{i+1} \not\equiv 0
\pmod{p}$ only for $i = p^{l} - 1$. 

Therefore for $0 < i < p^{l} - 1$ the summation~\eqref{eq:i}
vanishes. 

For $i = 0$ the summation~\eqref{eq:i} yields
\begin{math}
  - \dfrac{1}{2},
\end{math}
while for $i = p^{l} - 1$ it yields
\begin{equation*}
  (-1)^{i+1} \binom{\frac{\sigma + p^{l}}{2}}{i+1}
  =
  - \binom{p^{l} \frac{\sigma' + 1}{2}}{p^{l}}
  \equiv
  - \frac{\sigma' + 1}{2}
  \pmod{p}.
\end{equation*}
In conclusion, the summation in~\eqref{eq:summation} evaluates to
\begin{equation*}
  - \frac{1}{2} - \frac{\sigma' + 1}{2} = - \frac{\sigma' + 2}{2},
\end{equation*}
so that~\eqref{eq:summation} yields
\begin{equation*}
  0 
  = 
  [y \U{z}{\mu} y]
  - \frac{\sigma' + 2}{2}
  \,
  [y \U{z}{\mu + 1}]
  =
  - \frac{\sigma'}{2} 
  \,
  [y \U{z}{\mu} y],
\end{equation*}
as $[y \U{z}{\mu} x] = 0$. Since $\sigma' \not\equiv 0 \pmod{p}$, we
obtain a contradiction.

It is easy to see that one obtains exactly the same result when $p^{l}
= 1$.

\section{Specialization of two-step centralizers}
\label{sec:specialization}

In  this section  we will  prove Proposition~\ref{prop:specialization}
about specialization of centralizers.

So let $x_{1}$  and $x_{2}$ be as in  the statement of the Proposition. As
in  Section~3,  we  can clearly  write   all  two-step
centralizers other than $\F x_{1}$ in the form
\begin{equation}\label{eq:centr:F}
 \F \cdot ( x_{2} - \alpha x_{1} ),
\end{equation}
for some $\alpha \in \F$. Now let $\lambda$ be an indeterminate, and
consider the algebra
\begin{equation*}
 L ( \lambda ) = L \otimes_{\F} \F ( \lambda )
\end{equation*}
over  the  field  $\F  (  \lambda  )$.  We  can  redefine  $x_{1}$  as
$\lambda^{-1} x_{1}$  in $L  ( \lambda )$.   Thus $\F  (\lambda) \cdot
x_{1}$    and   $\F    (\lambda)   \cdot    x_{2}$    are   unchanged,
whereas~\eqref{eq:centr:F} becomes
\begin{equation}\label{eq:centr:Flambda}
 \F ( \lambda ) \cdot ( x_{2} - \alpha \lambda x_{1} ),
\end{equation}
with  $\alpha \in \F$.   As in  Section~3 of~I,  we can  describe $\ad
(x_{1})$ and $\ad (x_{2})$ by  introducing the following basis of $L (
\lambda )$:
\begin{equation*}
  v_{0} = x_{1},   \qquad
  v_{1} = x_{2}    \qquad
  v_{2} = [x_{2} x_{1}],
\end{equation*}
and for $i \ge 3$
\begin{equation}\label{eq:ad.a}
    v_{i} =
    \begin{cases}
      {}[v_{i-1} x_{1}] & \text{if $[v_{i-1} x_{1}] \ne 0$}\\
      {}[v_{i-1} x_{2}] & \text{if $[v_{i-1} x_{1}] = 0$.}
    \end{cases}
\end{equation}
When $[v_{i-1} x_{1}]$  is different
from zero, then $C_{L(\lambda)_{1}} ( v_{i-1} ) = \F ( \lambda ) \cdot
( x_{2} - \alpha \lambda x_{1} )$ for some $\alpha \in \F$, so that
\begin{equation}\label{eq:ad.b}
 [v_{i-1} x_{2}] = 
 [v_{i-1}, x_{2} - \alpha \lambda x_{1} + \alpha \lambda x_{1}] =
 \alpha \lambda [v_{i-1} x_{1}].
\end{equation}
It   follows  from~\eqref{eq:ad.a}   and   \eqref{eq:ad.b}  that   the
coefficients of the adjoint representation of $x_{1}$ and $x_{2}$ with
respect to the basis of the $v_{i}$ are in $\F [ \lambda ]$.  Consider
the  Lie subring  $S$ of  $L  ( \lambda  )$ generated  by $x_{1}$  and
$x_{2}$.  Then  the structure constants  of $S$ are in  the polynomial
ring $\F [ \lambda ]$. We can thus consider the algebra over $\F$
\begin{equation*}
  T = S \otimes_{\F [ \lambda ]} \trivfrac{\F [ \lambda ]}{( \lambda )}.
\end{equation*}
The relations~\eqref{eq:ad.a}  show that  $T$  is still an  algebra of
maximal class over $\F$,  which  is  effectively obtained  by  letting
$\lambda = 0$  in~\eqref{eq:centr:Flambda}.  The process of going from
$L$ to $T$ can thus be described in  terms of the respective sequences
of   \TSCs\ by saying  that we   are  leaving the  two-step centralizer
$x_{1}$ unchanged, and we are changing all other two-step centralizers
to $x_{2}$, as claimed.

\section{More than two distinct two-step centralizers}
\label{sec:more_than_two}
 
Using the result about specialization of two-step centralizers, we are
now    able    to   complete    the    classification,   by    proving
Step~\ref{step:more_than_two_is_inflated}, that is,  if the algebra of
maximal class  $L$ has three  or more distinct  two-step centralizers,
then   its   constituents   are    long   or   short,   so   that   by
Step~\ref{step:inflated} $L$ is inflated.

We first  prove Step~\ref{step:occurrence_of_centralizers}. Let $z$ be
any centralizer  higher than the first  two in order of occurrence. We
can leave $z$ fixed, and specialize all other  centralizers to $y$. In
the  resulting algebra,  $z$ will thus   play the role  of  the second
two-step   centralizer     in       order   of    occurrence.       By
Proposition~\ref{prop:parameter}, it will first occur as
\begin{equation*}
 C_{L_{1}} (L_{2 p^{n}}),
\end{equation*}
as stated.

Now we   prove  Lemma~\ref{lemma:third_after_short}.  Let $q$   be the
parameter of  $L$.  Again, let $z$  be any centralizer higher than the
first two in order of occurrence.  Suppose that  in the \SOTC\ we have
the segment
\begin{equation*}
 w_{1} \underbrace{y \dots y}_{m-1} z \underbrace{y \dots y}_{n-1} w_{2},
\end{equation*}
where $w_{1}, w_{2} \ne y$. Therefore $\underbrace{y \dots y}_{m-1} z$ and 
$\underbrace{y \dots y}_{n-1} w_{2}$ are two constituents, so that $m, n
\ge q$. Keep the second two-step centralizer $x$ fixed, and specialize all
others, including $z$, to $y$. This does not alter the parameter $q$. In the
resulting algebra, we have now a constituent that comprises at least
\begin{equation*}
 \underbrace{y \dots y}_{m-1} y \underbrace{y \dots y}_{n-1} w_{2}.
\end{equation*}
This has thus length at least $m +  n \ge 2 q$.  Since the length of a
constituent is at most $2 q$, it follows that $m =  n = q$, and $w_{1}
=      w_{2}      =       x$.       We      have        thus    proved
Lemma~\ref{lemma:third_after_short},   and    also  the other  related
results of~\cite{CaJur}.

We  are now able  to prove  Step~\ref{step:more_than_two_is_inflated}. 
Keep $x$ fixed in $L$,  and specialize all centralizers other than $x$
to $y$. As remarked, we are  not altering the parameter $q$. We obtain
an  algebra   of  maximal  class  $L'$  with   two  distinct  two-step
centralizers.  Because   of  Lemma~\ref{lemma:third_after_short},  the
first occurrence  of a centralizer higher  than the second  one in $L$
will  give   rise  to   a  long  constituent   in  $L'$.   Because  of
Step~\ref{step:pattern_for_AFS}, all the  constituents leading to this
occurrence in $L$ are either short or long. Hence the first non-short
constituent of $L'$ is long. It      follows      from
Step~\ref{step:Longs_and_Shorts} that  \emph{all} constituents in $L'$
are either short or long.  If we translate this from $L'$ back to $L$,
all that happens is that some of the long constituents of $L'$ are split
into two short  ones in $L$ by a two-step  centralizer higher than the
second  one. Thus  in $L$  too  all constituents  are short  or long.  
Step~\ref{step:more_than_two_is_inflated} is proved.

\providecommand{\bysame}{\leavevmode\hbox to3em{\hrulefill}\thinspace}


\begin{thebibliography}{CMN97}

\bibitem[Car98]{claretta}
C.~Carrara, \emph{({F}inite) presentations of loop algebras of {Albert-Frank
  Lie} algebras}, Ph.D. thesis, Trento, 1998.

\bibitem[Car99]{claretta2}
C.~Carrara, \emph{({F}inite) presentations of loop algebras of {Albert-Frank
  Lie} algebras}, submitted, 1999.

\bibitem[CJ99]{CaJur}
A.~Caranti and G.~Jurman, \emph{Quotients of maximal class of thin {Lie}
  algebras. {T}he odd characteristic case}, accepted for publication in Comm.\
  Algebra, 1999.

\bibitem[CMN97]{max}
A.~Caranti, S.~Mattarei, and M.~F. Newman, \emph{Graded {L}ie algebras of
  maximal class}, Trans. Amer. Math. Soc. \textbf{349} (1997), no.~10,
  4021--4051.

\bibitem[HNO97]{HNO}
George Havas, M.F. Newman, and E.A. O'Brien, \emph{{ANU} $p$-{Q}uotient
  {P}rogram (version 1.4)}, written in \texttt{C}, available as a share library
  with \textsf{GAP} and as part of \texttt{{M}agma}, or from
  {\texttt{http://wwwmaths.anu.edu.au/services/ftp.html}}, School of
  Mathematical Sciences, Australian National University, Canberra, 1997.

\bibitem[Jur98]{beppe2}
G.~Jurman, \emph{On graded {Lie} algebras in characteristic two}, Ph.D. thesis,
  Trento, 1998.

\bibitem[Jur99]{beppe}
G.~Jurman, \emph{Graded {Lie} algebra of maximal class. {III}}, in preparation,
  1999.

\bibitem[KW89]{Knuth:Lucas}
Donald~E. Knuth and Herbert~S. Wilf, \emph{The power of a prime that divides a
  generalized binomial coefficient}, J. reine angew. Math. \textbf{396} (1989),
  212--219.

\bibitem[Luc78]{Lucas}
\`E. Lucas, \emph{Sur les congruences des nombres eul\'eriens et des
  coefficients diff\'erentiels des fonctions trigonom\'etriques, suivant un
  module premier}, Bull. Soc. Math. France \textbf{6} (1878), 49--54.

\bibitem[Sha94]{Sha:max}
Aner Shalev, \emph{Simple {L}ie algebras and {L}ie algebras of maximal class},
  Arch. Math. (Basel) \textbf{63} (1994), no.~4, 297--301.

\end{thebibliography}
\end{document}